\documentclass{article}%
\usepackage{amsmath}
\usepackage{graphicx}
\usepackage{amsfonts}
\usepackage{amssymb}%
\newtheorem{theorem}{Theorem}

\newtheorem{definition}[theorem]{Definition}
\newtheorem{example}[theorem]{Example}

\newtheorem{proposition}[theorem]{Proposition}

\newenvironment{proof}[1][Proof]{\textbf{#1.} }{\ \rule{0.5em}{0.5em}}
\begin{document}

\title{Geodesic Webs and PDE Systems of Euler Equations}
\author{Vladislav V. Goldberg\\New Jersey Institute of Technology, Newark, NJ, USA
\and Valentin V. Lychagin\\University of Tromso, Tromso, Norway}
\date{}
\maketitle

\begin{abstract}
We find necessary and sufficient conditions for the foliation defined by level
sets of a function $f(x_{1},...,x_{n})$ to be totally geodesic in a
torsion-free connection and apply them to find the conditions for $d$-webs of
hypersurfaces to be geodesic, and in the case of flat connections, for
$d$-webs ($d\geq n+1$) of hypersurfaces to be hyperplanar webs. These
conditions are systems of generalized Euler equations, and for flat
connections we give an explicit construction of their solutions.

\end{abstract}

\section{Introduction}

In this paper we study  necessary and sufficient conditions for the foliation
defined by level sets of a function to be totally geodesic in a torsion-free
connection on a manifold and find necessary and sufficient conditions for webs
of hypersurfaces to be geodesic. These conditions has the form of a second-order
PDE system for web functions. The system has an infinite
pseudogroup of symmetries and the factorization of the system with respect to
the pseudogroup leads us to a first-order PDE system. In the planar
case (cf. \cite{GL2008}), the system coincides with the classical Euler
equation and therefore can be solved in a constructive way. We provide a
method to solve the system in arbitrary dimension and flat connection.

\section{Geodesic Foliations and Flex Equations}

Let $M^{n}$ be a smooth manifold of dimension $n.$ Let vector fields
$\partial_{1},...,\partial_{n}$ form a basis in the tangent bundle, and let
$\omega^{1},..,\omega^{n}$ be the dual basis. Then
\[
\lbrack\partial_{i},\partial_{j}]=\sum_{k}c_{ij}^{k}\partial_{k}%
\]
for some functions $c_{ij}^{k}\in C^{\infty}\left(  M\right)  ,$ and
\[
d\omega^{k}+\sum_{i<j}c_{ij}^{k}\omega^{i}\wedge\omega^{j}=0.
\]
Let $\nabla$ be a linear connection in the tangent bundle, and let
$\Gamma_{ij}^{k}$ be the Christoffel symbols of second type. Then%
\[
\nabla_{i}\left(  \partial_{j}\right)  =\sum_{k}\Gamma_{ij}^{k}\partial_{k},
\]
where $\nabla_{i}\overset{\text{def}}{=}\nabla_{\partial_{i}},$ and
\[
\nabla_{i}\left(  \omega^{k}\right)  =-\sum_{j}\Gamma_{ij}^{k}\omega^{j}.
\]

In \cite{GL2008} we proved the following result.

\begin{theorem}
The foliation defined by the level sets of a function $f(x_{1},...,x_{n})$ is
totally geodesic in a torsion-free connection $\nabla$ if and only if the
function $f$ satisfies the following system of PDEs:%
\begin{equation}
\frac{\partial_{i}\left(  f_{i}\right)  }{f_{i}f_{i}}-\frac{\partial
_{i}\left(  f_{j}\right)  +\partial_{j}\left(  f_{i}\right)  }{f_{i}f_{j}%
}+\frac{\partial_{j}\left(  f_{j}\right)  }{f_{j}f_{j}}=\sum_{k}\left(
\Gamma_{ii}^{k}\frac{f_{k}}{f_{i}f_{i}}+\Gamma_{jj}^{k}\frac{f_{k}}{f_{j}%
f_{j}}-(\Gamma_{ij}^{k}+\Gamma_{ji}^{k})\frac{f_{k}}{f_{i}f_{j}}\right)
\label{geodcondfor f}%
\end{equation}
for all $i<j,i,j=1,...,n;$ here $f_{i}=\displaystyle           \frac{\partial
f}{\partial x_{i}}.$
\end{theorem}

We call such a system a \emph{flex system.}

Note that conditions (\ref{geodcondfor f}) can be used to obtain necessary and
sufficient conditions for a $d$-web formed by the level sets of the functions
$f_{\alpha}(x_{1},\dots,x_{n}),\alpha=1,\dots,d,$ to be a \emph{geodesic }%
$d$-\emph{web,} i.e., to have the leaves of all its foliations to be totally
geodesic: one should apply conditions (\ref{geodcondfor f}) to the all web
functions $f_{\alpha},\alpha=1,\dots,d,$

\subsection{Geodesic Webs on Manifolds of Constant Curvature}

In what follows, we shall use the following definition.

\begin{definition}
We call by $(\operatorname{Flex}\;f)_{ij}$ the following function:%
\[
(\operatorname{Flex}\;f)_{ij}=f_{j}^{2}f_{ii}-2f_{i}f_{j}f_{ij}+f_{i}%
^{2}f_{jj},
\]
where $i,j=1,...,n,$ $f_{i}=\displaystyle\frac{\partial f}{\partial x_{i}}$
and $f_{ij}=\displaystyle\frac{\partial^{2}f}{\partial x_{i}\partial x_{j}}.$
\end{definition}

It is easy to see that $(\operatorname{Flex}\;f)_{ij}=(\operatorname{Flex}%
\;f)_{ji},$ and $(\operatorname{Flex}\;f)_{ii}=0.$

\begin{proposition}
Let $\left(  \mathbb{R}^{n},g\right)  $ be a manifold of constant curvature
with the metric tensor
\[
g=\frac{dx_{1}^{2}+...+dx_{n}^{2}}{\left(  1+\kappa\left(  x_{1}^{2}%
+...+x_{n}^{2}\right)  \right)  ^{2}},
\]
where $\kappa$ is a constant. Then the level sets of a function $f(x_{1}%
,...,x_{n})$ are geodesics of the metric $g$ if and only if the function $f$
satisfies the following PDE system:
\begin{equation}
(\operatorname{Flex}\;f)_{ij}=\frac{2\kappa\left(  f_{i}^{2}+f_{j}^{2}\right)
}{1+\kappa\left(  x_{1}^{2}+...+x_{n}^{2}\right)  }\sum_{k}x_{k}%
f_{k}\label{flex equation}%
\end{equation}
for all $i,j.$
\end{proposition}

\begin{proof}
To prove formula (\ref{flex equation}), first note that the components of the
metric tensor $g$ are
\[
g_{ii}=b^{2},\;g_{ij}=0,\;\;i\neq j,
\]
where
\[
b=\frac{1}{1+\kappa\left(  x_{1}^{2}+...+x_{n}^{2}\right)  }.
\]
It follows that%
\[
g^{ii}=g_{ii}^{-1},\;g^{ij}=0,\;\;i\neq j.
\]

We compute $\Gamma_{jk}^{i}$ using the classical formula%
\begin{equation}
\renewcommand{\arraystretch}{1.3}\Gamma_{ij}^{k}=\frac{1}{2}g^{kl}\left(
\frac{\partial g_{li}}{\partial x^{j}}+\frac{\partial g_{lj}}{\partial x^{i}%
}-\frac{\partial g_{ij}}{\partial x^{l}}\right)
\renewcommand{\arraystretch}{1} \label{Christoffel}%
\end{equation}
and get%
\[
\renewcommand{\arraystretch}{1.3}%
\begin{array}
[c]{ll}%
\Gamma_{ii}^{k}=2\kappa x_{k}b,\;k\neq i;\;\Gamma_{ii}^{i}=-2\kappa
x_{i}b;\;\Gamma_{ij}^{k}=0,\;i,j\neq k,\;i\neq j;\; & \\
\Gamma_{ij}^{i}=-2\kappa x_{j}b,\;i\neq j;\;\Gamma_{ij}^{j}=-2\kappa
x_{i}b,\;\;i\neq j. &
\end{array}
\renewcommand{\arraystretch}{1}
\]

Substituting these values of $\Gamma_{jk}^{i}$ into the right-hand side of
formula (\ref{geodcondfor f}), we get formula (\ref{flex equation}).
\end{proof}

Note that if $n=2$, then PDE system (\ref{flex equation}) reduces to the
single equation
\[
\operatorname{Flex}\;f=\frac{2\kappa\left(  x_{1}f_{1}+x_{2}f_{2}\right)
\left(  f_{1}^{2}+f_{2}^{2}\right)  }{1+\kappa\left(  x_{1}^{2}+x_{2}%
^{2}\right)  },
\]
where $\operatorname{Flex}\;f=(\operatorname{Flex}\;f)_{12}$.

This formula coincides with the corresponding formula in \cite{GL2008}.

We rewrite formula (\ref{flex equation}) as follows:%

\begin{equation}
\displaystyle\frac{(\operatorname{Flex}\;f)_{ij}}{f_{i}^{2}+f_{j}^{2}}=2\kappa
b\sum_{k}x_{k}f_{k}. \label{flex equation2}%
\end{equation}
The left-hand side of equation (\ref{flex equation2}) does not depend on $i$
and $j.$ Thus we have%

\[
\displaystyle\frac{(\operatorname{Flex}\;f)_{ij}}{f_{i}^{2}+f_{j}^{2}%
}=\displaystyle\frac{(\operatorname{Flex}\;f)_{kl}}{f_{k}^{2}+f_{l}^{2}}%
\]
for any $i,j,k,$ and $l.$

It follows that \textit{if }%
\begin{equation}
(\operatorname{Flex}\;f)_{ij}=0\label{flex=0}%
\end{equation}
\textit{for some fixed }$i\;$\textit{and}$\;j,\;$\textit{then} (\ref{flex=0})
\textit{holds for any }$i$ \textit{and}$\;j.\;$

In other words, one has the following result.

\begin{theorem} \label{planar intersections}
Let $W$ be a geodesic $d$-web on the manifold $\left(  \mathbb{R}%
^{n},g\right)  $ given by web-functions $\left\{  f^{1},....,f^{d}\right\}  $
such that $\left(  f_{k}^{a}\right)  ^{2}+\left(  f_{l}^{a}\right)  ^{2}\neq 0$
for all $a=1,...,d$ and $k,l=1,2...,n.$ Assume that the intersections of $W$ with the
planes $\left(  x_{i_{0}},x_{j_{0}}\right),$ for given $i_{0}$ and $j_{0},$ are
linear planar $d$-webs. Then the intersection of $W$ with arbitrary planes $\left(
x_{i},x_{j}\right)  $ are linear webs too.
\end{theorem}

\subsection{Geodesic Webs on Hypersurfaces in $\mathbb{R}^{n}$}

\begin{proposition}
Let $\left(  M,g\right)  \subset\mathbb{R}^{n}$ be a hypersurface defined by
an equation $x_{n}=u\left(  x_{1},...,x_{n-1}\right)  $ with the induced
metric $g$ and the Levi-Civita connection $\nabla$. Then the foliation defined
by the level sets of a function $f\left(  x_{1},...,x_{n-1}\right)  $ is
totally geodesic in the connection $\nabla$ if and only if the function $f$
satisfies the following system of PDEs:
\begin{equation}
(\operatorname{Flex}\;f)_{ij}=\frac{u_{1}f_{1}+...+u_{n-1}f_{n-1}}{1+u_{1}%
^{2}+...+u_{n-1}^{2}}(f_{j}^{2}u_{ii}-2f_{i}f_{j}u_{ij}+f_{i}^{2}u_{jj}).
\label{flex for u(x1,...,xn-1)}%
\end{equation}
\end{proposition}

\begin{proof}
To prove formula (\ref{flex for u(x1,...,xn-1)}), note that the metric induced
by a surface $x_{n}=u(x_{1},\dots,x_{n-1})$ is
\[
g=ds^{2}=\sum_{k=1}^{n-1}(1+u_{k}^{2})dx_{k}^{2}+2\sum_{i,j=1(i\neq j)}%
^{n-1}u_{i}u_{j}dx_{i}dx_{j}.
\]
Thus the metric tensor $g$ has the following matrix:%
\[
\renewcommand{\arraystretch}{1.5} (g_{ij})=\left(
\begin{array}
[c]{cccccccccccccc}%
\displaystyle 1+u_{1}^{2} & u_{1}u_{2} & \dots & u_{1}u_{n-1}\cr u_{2}u_{1} &
1+u_{2}^{2} & \dots & u_{2}u_{n-1}\cr \vdots & \vdots & \ddots &
\vdots\cr u_{1} & u_{n-1}u_{2} & \dots & 1+u_{n-1}^{2}, &
\end{array}
\right)  \renewcommand{\arraystretch}{1}
\]
and the inverse tensor $g^{-1}$ has the matrix
\[
(g^{ij})=\displaystyle        \frac{1}{1+\displaystyle   \sum_{k=1}%
^{n-1}(1+u_{k}^{2})}\left(
\begin{array}
[c]{cccccccccccccc}%
\displaystyle \sum_{k=2}^{n-1}(1+u_{k}^{2}) & -u_{1}u_{2} & \dots &
-u_{1}u_{n-1}\cr -u_{2}u_{1} & \displaystyle \sum_{k=1(k\neq2)}^{n-1}%
(1+u_{k}^{2}) & \dots & -u_{2}u_{n-1}\cr \vdots & \vdots & \ddots &
\vdots\cr -u_{n-1}u_{1} & -u_{n-1}u_{2} & \dots & \displaystyle \sum
_{k=1}^{n-2}(1+u_{k}^{2}) &
\end{array}
\right)  .\renewcommand{\arraystretch}{1}
\]
Computing $\Gamma_{jk}^{i}$ by formula (\ref{Christoffel}), we find that%
\[
\Gamma_{ij}^{k}=\displaystyle               \frac{u_{k}u_{ij}}%
{1+\displaystyle                                     \sum_{k=1}^{n-1}%
(1+u_{k}^{2})}.
\]
Applying these formulas to the right-hand side of (\ref{geodcondfor f}), we
get formula (\ref{flex for u(x1,...,xn-1)}).
\end{proof}

We rewrite equation (\ref{flex for u(x1,...,xn-1)}) in the form
\begin{equation}
\displaystyle \frac{(\operatorname{Flex}\;f)_{ij}}{f_{j}^{2}u_{ii}-2f_{i}%
f_{j}u_{ij}+f_{i}^{2}u_{jj}}=\frac{u_{1}f_{1}+...+u_{n}f_{n}}{1+u_{1}%
^{2}+...+u_{n}^{2}} . \label{flex2 for u(x1,...,xn-1)}%
\end{equation}

It follows that the left-hand side of (\ref{flex2 for u(x1,...,xn-1)}) does
not depend on $i$ and $j$, i.e., we have
\[
\displaystyle\frac{(\operatorname{Flex}\;f)_{ij}}{f_{j}^{2}u_{ii}-2f_{i}%
f_{j}u_{ij}+f_{i}^{2}u_{jj}}=\displaystyle\frac{(\operatorname{Flex}\;f)_{kl}%
}{f_{l}^{2}u_{kk}-2f_{k}f_{l}u_{kl}+f_{k}^{2}u_{ll}}%
\]
for any $i,j,k$ and $l$. This means that \textit{if
\[
(\operatorname{Flex}\;f)_{ij}=0~
\]
for some fixed $i$ and $j$, then
\[
(\operatorname{Flex}\;f)_{kl}=0
\]
for any $k$ and $l$}.

In other words, we have a result similar to the result in Theorem \ref{planar intersections}.

\begin{theorem}
Let $W$ be a geodesic $d$-web on the hypersurface  $\left(  M,g\right)  $
given by web functions $\left\{  f^{1},....,f^{d}\right\}  $ such that
$\left(  f_{j}^{a}\right)  u_{ii}-2f_{i}^{a}f_{j}^{a}u_{ij}+\left(  f_{i}%
^{a}\right)^{2}u_{jj}\neq 0$, for all $a=1,...,d$ and $k,l=1,2...,n.$ Assume
that the intersections of $W$ with the planes $\left(  x_{i_{0}},x_{j_{0}}\right)  ,$ for
given $i_{0}$ and $j_{0},$ are linear planar $d$-webs. Then the intersection of $W$ with
arbitrary planes $\left(  x_{i},x_{j}\right)  $ are linear webs too.
\end{theorem}

\section{Hyperplanar Webs}

In this section we consider hyperplanar geodesic webs in \ $\mathbb{R}^{n}$
endowed with a flat linear connection $\nabla$.

In what follows, we shall use coordinates $x_{1},\dots,x_{n}$ in which the
Christoffel symbols $\Gamma_{jk}^{i}$ of $\nabla$ vanish.

The following theorem gives us a criterion for a web of hypersurfaces to be hyperplanar.

\begin{theorem}
Suppose that a $d$-web of hypersurfaces, $d\geq n+1$, is given locally by web
functions $f_{\alpha}(x_{1},\dots,x_{n}),\alpha=1,...,d.$ Then the web is
hyperplanar if and only if the web functions satisfy the following PDE
system:
\begin{equation}
(\operatorname{Flex}\;f)_{st}=0, \label{linear d-web}%
\end{equation}
for all $\;s<t=1,...,n.$
\end{theorem}

\begin{proof}
For the proof, one should apply Theorem 1 to all foliations of the web.
\end{proof}

In order to integrate the above PDEs system, we introduce the functions%
\[
A_{s}{}=\frac{f_{s}}{f_{s+1}},\;\;s=1,...,n-1,
\]
and the vector fields%
\[
X_{s}=\frac{\partial}{\partial x_{s}}-A_{s}\frac{\partial}{\partial x_{s+1}%
},\;s=1,...,n-1.
\]
Then the system can be written as
\[
X_{s}\left(  A_{t}\right)  =0,
\]
where $s,t=1,..,n-1.$

Note that
\[
\lbrack X_{s},X_{t}]=0
\]
if the function $f$ is a solution of (\ref{linear d-web}).

Hence, the vector fields $X_{1},...,X_{n-1}$ generate a completely integrable
$\left(  n-1\right)  $-dimensional distribution, and the functions
$A_{1},...,A_{n-1}$ are the first integrals of this distribution.

Moreover, the definition of the functions $A_{s}$ shows that
\[
X_{s}(f)=0,\ s=1,...,n-1,
\]
also.

As a result, we get that
\[
A_{s}=\Phi_{s}\left(  f\right)  ,\ \ s=1,...,n-1,
\]
for some functions $\Phi_{s}.$

In these terms, we get the following system of equations for $f$:%
\[
\frac{\partial f}{\partial x_{s}}=\Phi_{s}\left(  f\right)  \frac{\partial
f}{\partial x_{s+1}},\,\ s=1,...,n-1,
\]
or%
\begin{equation}%
\begin{array}
[c]{l}%
\displaystyle\frac{\partial f}{\partial x_{s}}=\Psi_{s}\left(  f\right)
\displaystyle\frac{\partial f}{\partial x_{n}},\;s=1,...,n-1,
\end{array}
\label{system for webfunction f of linear (n+1)-web)}%
\end{equation}
where $\Psi_{n-1}=\Phi_{n-1},$ and
\[
\Psi_{s}=\Phi_{n-1}\cdots\Phi_{s}%
\]
for $s=1,....,n-2.$

This system is a sequence of the Euler-type equations and therefore can be
integrated. Keeping in mind that a solution of the single Euler-type equation%
\[
\frac{\partial f}{\partial x_{s}}=\Psi_{s}\left(  f\right)  \frac{\partial
f}{\partial x_{n}}%
\]
is given by the implicit equation
\[
f=u_{0}\left(  x_{n}+\Psi_{s}\left(  f\right)  x_{s}\right)  ,
\]
where $u_{0}(x_{n})$ is an initial condition, when $x_{s}=0,\;$and $\Psi_{s}$
is an arbitrary nonvanishing function, we get solutions $f$ of system
(\ref{linear d-web}) in the form:
\[
f=u_{0}\left(  x_{n}+\Psi_{n-1}\left(  f\right)  x_{n-1}+\cdots+\Psi
_{1}\left(  f\right)  x_{1}\right)  ,
\]
where $u_{0}(x_{n})$ is an initial condition, when $x_{1}=\cdots=x_{n-1}%
=0,\;$and $\Psi_{s}$ are arbitrary nonvanishing functions.

Thus, we have proved the following result.

\begin{theorem}
Web functions of hyperplanar webs have the form
\begin{equation}
f=u_{0}\left(  x_{n}+\Psi_{n-1}\left(  f\right)  x_{n-1}+\cdots+\Psi
_{1}\left(  f\right)  x_{1}\right)  ,
\label{webfunction f for linear (n+1)-web}%
\end{equation}
where $u_{0}(x_{n})$ are initial conditions, when $x_{1}=\cdots=x_{n-1}%
=0,\;$and $\Psi_{s}$ are arbitrary nonvanishing functions.
\end{theorem}

\begin{example}
\label{excone} Assume that $n=3,\;f_{1}(x_{1},x_{2},x_{3})=x_{1},\;f_{2}%
(x_{1},x_{2},x_{3})=x_{2},$ $f_{3}(x_{1},x_{2},x_{3})=x_{3}$, and take
$u_{0}=x_{3},\;\Psi_{1}(f_{4})=f_{4}^{2},\;\Psi_{2}(f_{4})=f_{4}$ in
$(\ref{webfunction f for linear (n+1)-web})$. Then we get the hyperplanar
$4$-web with the remaining web function
\[
f_{4}=\displaystyle\frac{x_{2}-1\pm\sqrt{(x_{2}-1)^{2}-4x_{1}x_{3}}}{2x_{1}}.
\]
It follows that the level surfaces $f_{4}=C$ of this function are defined by
the equation
\[
x_{1}(C^{2}x_{1}-Cx_{2}+x_{3}+C)=0,
\]
i.e., they form a one-parameter family of $2$-planes%
\[
C^{2}x_{1}-Cx_{2}+x_{3}+C=0.
\]
Differentiating the last equation with respect to $C$ and excluding $C,$ we
find that the envelope of this family is defined by the equation%
\[
(x_{2})^{2}-4x_{1}x_{3}-2x_{2}+1=0.
\]
Therefore, the envelope is the second-degree cone.
\end{example}

\begin{example}
\label{exhypcylinder} Assume that $n=3,\;f_{1}(x_{1},x_{2},x_{3}%
)=x_{1},\;f_{2}(x_{1},x_{2},x_{3})=x_{2},$ $f_{3}(x_{1},x_{2},x_{3})=x_{3}$,
and take $u_{0}=x_{3},\Psi_{1}(f_{4})=1,\Psi_{2}(f_{4})=f_{4}^{2}$ in
$(\ref{webfunction f for linear (n+1)-web})$. Then we get the linear $\,4$-web
with the remaining web function
\[
f_{4}=\Biggl(\displaystyle\frac{1\pm\sqrt{1-4x_{2}(x_{1}+x_{3})}}{2x_{2}%
}\Biggr)^{2}.
\]
The level surfaces $f_{4}=C^{2}$ of this function are defined by the equation
\[
x_{2}(x_{1}+C^{2}x_{2}+x_{3}-C)=0,
\]
i.e., they form a one-parameter family of $2$-planes%
\[
x_{1}+C^{2}x_{2}+x_{3}-C=0.
\]
Differentiating the last equation with respect to $C$ and excluding $C,$ we
find that the envelope of this family is defined by the equation%
\[
4x_{1}x_{2}+4x_{2}x_{3}-1=0.
\]
Therefore, the envelope is the hyperbolic cylinder.
\end{example}

In the next example no one foliation of a web $W_{3}\;$coincides with a
foliation of coordinate lines, i.e., all three web functions are unknown.

\begin{example}
Assume that $n=3$ and take

\begin{description}
\item[(i)]
$u_{01}=x_{3},\;\Psi_{1}(f_{1})=f_{1}^{2},\;\Psi_{2}(f_{1})=f_{1}$;

\item[(ii)]
$u_{02}=x_{3},\Psi_{1}(f_{2})=1,\Psi_{2}(f_{2})=f_{2}^{2};$

\item[(iii)]
$u_{03}=x_{3}^{2},\;\Psi_{1}(f_{3})=f_{3},\;\Psi_{2}(f_{3})=1;$

\item[(iv)] $u_{04}=x_{3},\;\Psi_{1}(f_{4})=\Psi_{2}(f_{4})=f_{4}$
\end{description}

\noindent in $(\ref{webfunction f for linear (n+1)-web})$. Then we
get the linear $4$-web with the web functions
\begin{align*}
f_{1}  &  =\displaystyle\frac{x_{2}-1\pm\sqrt{(x_{2}-1)^{2}-4x_{1}x_{3}}%
}{2x_{1}},\\
f_{2}  &  =\Biggl(\displaystyle\frac{1\pm\sqrt{1-4x_{2}(x_{1}+x_{3})}}{2x_{2}%
}\Biggr)^{2}%
\end{align*}
$($see Examples $\ref{excone}$ and $\ref{exhypcylinder})$ and
\begin{align*}
f_{3}  &  =(\frac{1\pm\sqrt{1-4x_{1}(x_{2}+x_{3})}}{2x_{1}})^{2},\\
f_{4}  &  =\frac{x_{3}}{1-x_{1}-x_{2}}.
\end{align*}

It follows that the leaves of the foliation $X_{1}$ are tangent
$2$-planes to the second-degree cone
\[
(x_{2})^{2}-4x_{1}x_{3}-2x_{2}+1=0
\]
$($cf. Example $\ref{excone}$ and $\ref{exhypcylinder}),$ the
leaves of the foliation $X_{2}$ and $X_{3}$ are tangent $2$-planes
to the hyperbolic cylinders
\[
4x_{1}x_{2}+4x_{2}x_{3}-1=0 \;\text{and} \;4x_{1}x_{2}+4x_{1}x_{3}-1=0
\]
$($cf. Example $\ref{exhypcylinder}),$ and the leaves of the
foliation $X_{4}$
are $2$-planes of the one-parameter family of parallel $2$-planes%
\[
Cx_{1}+Cx_{2}+x_{3}=1,
\]
where $C$ is an arbitrary constant.
\end{example}

{\emph{Authors' addresses:} }

{Deparment of Mathematical Sciences, New Jersey Institute of Technology,
University Heights, Newark, NJ 07102, USA; vladislav.goldberg@gmail.com}

{Department of Mathematics, The University of Tromso, N9037,
Tromso, Norway; lychagin@math.uit.no}

\end{document}